\documentclass{amsart}
\usepackage{amsmath}

\newtheorem{theorem}{\textbf{Theorem}}
\newtheorem{question}{\textbf{Question}}

\def\A {\mathbb{A}}

\def\a {\alpha}

\def\Q {\mathbb{Q}}
\def\R {\mathbb{R}}

\def\QQ {\overline{\Q}}
\def\C {\mathbb{C}}

\def\e {\epsilon}

\theoremstyle{remark}

\numberwithin{equation}{section}

\begin{document}

\title[A POSITIVE ANSWER FOR A QUESTION PROPOSED BY K. MAHLER]{A POSITIVE ANSWER FOR A QUESTION PROPOSED BY K. MAHLER}

\author{DIEGO MARQUES}
\address{DEPARTAMENTO DE MATEM\'{A}TICA, UNIVERSIDADE DE BRAS\'ILIA, BRAS\'ILIA, DF, BRAZIL}
\email{diego@mat.unb.br}


\author{CARLOS GUSTAVO MOREIRA}
\address{INSTITUTO DE MATEM\' ATICA PURA E APLICADA, RIO DE JANEIRO, RJ, BRAZIL}
\email{gugu@impa.br}



\begin{abstract}
In 1902, P. St\"{a}ckel proved the existence of a transcendental function $f(z)$, analytic in a neighbourhood of the origin, and with the property that both $f(z)$ and its inverse function assume, in this neighbourhood, algebraic values at all algebraic points. Based on this result, in 1976, K. Mahler raised the question of the existence of such functions which are analytic in $\mathbb{C}$. In this work, we provide a positive answer for this question by showing the existence of uncountable many of these functions.
\end{abstract}

\maketitle

\section{Introduction}

A \textit{transcendental function} is a function $f(x)$ such that the only complex polynomial satisfying $P(x, f(x)) =0$ for all $x$ in its domain, is the null polynomial. For instance, the trigonometric functions, the exponential function, and their inverses. 

The study of the arithmetic behavior of transcendental functions at complex points has attracted the attention of many mathematicians for decades. The first result concerning this subject goes back to 1884, when Lindemann proved that the transcendental function $e^z$ assumes transcendental values at all nonzero algebraic point. In 1886, Strauss tried to prove that an analytic transcendental function  cannot assume rational values at all rational points in its domain. However, in 1886, Weierstrass supplied him with a counter-example and also stated that there are transcendental entire functions which assume algebraic values at all algebraic points. This assertion was proved in 1895 by St\"{a}ckel \cite{19} who established a much more general result: for each countable subset $\Sigma\subseteq \C$ and each dense subset $T\subseteq \C$, there exists a transcendental entire
function $f$ such that $f(\Sigma) \subseteq T$. In another construction, St\"{a}ckel \cite{20} produced a transcendental function $f(z)$, analytic in a neighbourhood of the origin, and with the property that both $f(z)$ and its inverse function assume, in this neighbourhood, algebraic values at all algebraic points. Based on this result, in 1976, Mahler \cite[p. 53]{bookmahler} suggested the following question

\begin{question}\label{q1}
Does there exist a transcendental entire function
\[
f(z)=\sum_{n=0}^{\infty}a_nz^n,
\]
with rational coefficients $a_n$ and such that the image and the preimage of $\QQ$ under $f$ are subsets of $\QQ$?
\end{question}

We refer the reader to \cite{bookmahler,wal1} (and references therein) for more about this subject.

In this paper, we give an affirmative answer to Mahler's question. For the sake of preciseness, we stated it as a theorem.
\begin{theorem}\label{1}
There are uncountable many transcendental entire functions
\[
f(z)=\sum_{n=0}^{\infty}a_nz^n,
\]
with rational coefficients $a_n$ and such that the $f(\QQ)\subseteq \QQ$ and $f^{-1}(\QQ)\subseteq \QQ$.
\end{theorem}

\section{The proof}

In order to simplify our presentation, we set $\A=\QQ\cap \R$ and we use the familiar notation $[a, b] = \{a, a + 1,\ldots, b\}$, for integers $a < b$.

Let $\{\a_1,\a_2,\a_3,\ldots\}$ be an enumeration of $\QQ$ such that $\a_1=0$, and for any $n\geq 1$, the numbers $\a_{3n-1},\a_{3n}\notin \R$ with $\a_{3n}=\overline{\a}_{3n-1}$, $\a_{3n+1}\in \R$ and $|\a_{3n+i}|<n$, for $i\in [-1,1]$. Now, let us construct our desired function inductively. 

Define $f_1(z)=z$. Observe that $f_1(\a_1)=0=(f_1)^{-1}(\a_1)$. Now, we want to construct a sequence of analytic functions $f_2(z), f_3(z),\ldots$ recursively of the form
\[
f_m(z)=f_{m-1}(z)+\e_mz^mP_m(z),
\]
with
\begin{itemize}
\item[(i)] $f_{m-1},P_m\in \A[z]$ and then $f_m(z)=\sum_{i=1}^{t_m}a_iz^i$ with $t_m>m$;
\item[(ii)] $P_{m-1}(z)\mid P_m(z)$ and $P_m(0)\neq 0$;
\item[(iii)] $\e_m\in \A$;
\item[(iv)] $0<|\e_m|<\frac{1}{L(P_m)m^{m+\deg P_m}}$;
\item[(v)] $a_1,\ldots, a_m\in \Q$;
\end{itemize}
Here $L(P)$ denotes the length of a polynomial (the sum of the absolute values of its coefficients).

The requested function will have the form $f(z)=z+\sum_{n\geq 2}\e_nz^nP_n(z)$. Since $|P(z)|\leq L(P)\max\{1,|z|\}^{\deg P}$, we have that, for all $z$ belonging to the open ball $B(0,R)$,
\[
|\e_nz^nP_n(z)|<\dfrac{1}{L(P_n)n^{n+\deg P_n}}L(P_n)\max\{1,R\}^{n+\deg P_n}=\left(\dfrac{\max\{1,R\}}{n}\right)^{n+\deg P_n}.
\]
Thus $f$ is an entire function, since the series $f(z)=z+\sum_{n=2}^{\infty}\e_nz^nP_n(z)$, which defines $f$, converges uniformly in any of these balls.

Suppose that we have a function $f_{n}$ satisfying (i)-(v). Now, let us construct $f_{n+1}$ with the desired properties.

Since $f_{n}^{-1}(\{\a_1,\ldots, \a_{3n+1}\})$ is finite, we can choose $n+1<r_{n+1}<n+2$ such that $f_{n}^{-1}(\{\a_1,\ldots, \a_{3n+1}\})\cap \partial B(0:r_{n+1})=\emptyset$. If $f_{n}^{-1}(\{\a_1,\ldots, \a_{3n+1}\})\cap B(0:r_{n+1})=\{0,y_1,\ldots, y_s\}$. Define
\[
f_{n+1}(z)=f_{n}(z)+\e_{n+1}z^{n+1}P_{n+1}(z),
\]
where $P_{n+1}(z):=P_n(z)(z-\a_{3n-1})(z-\a_{3n})(z-\a_{3n+1})\prod_{i=1}^s(z-y_i)^{\deg f_n}$. Note that since $f_n(z)\in \A[z]$, one has that if $y_i\notin \R$, then $\overline{y_i}=y_j$, for some $j\in [1,s]$ (since $f_{n+1}(y_i)=f_n(y_i)$ and $\{\a_1,\ldots, \a_{3n+1}\}$ is closed under complex conjugation). Thus $P_{n+1}\in \A[z]$ and $P_{n+1}(0)\neq 0$.


For $i\in [1,3n+1]$, since $f_{n}^{-1}(\a_i)\cap \partial B(0:r_{n+1})$ is empty, then $\min_{|z|=r_{n+1}} |f_{n}(z)-\a_i|>0$ and we can choose $\e_{n+1}$ satisfying
\begin{equation}\label{rouche0}
|\e_{n+1}|<\dfrac{\displaystyle\min_{|z|=r_{n+1}}|f_n(z)-\a_i|}{\displaystyle\max_{|z|=r_{n+1}}|z^{n+1}P_{n+1}(z)|}.
\end{equation}
In particular,
\begin{eqnarray}\label{rouche}
|f_{n}(z)-\a_i| & \geq & \min_{|z|=r_{n+1}} |f_{n}(z)-\a_i|\nonumber \\
 & >& |\e_{n+1}|\max_{|z|=r_{n+1}}|z^{n+1}P_{n+1}(z)|\nonumber \\
 & \geq &  |\e_{n+1}z^{n+1}P_{n+1}(z)|.\nonumber 
\end{eqnarray}
Now, by using Rouch\' e's theorem, we ensure that the functions $f_n(z)-\a_i$ and $f_n(z) -\a_i+ \e_{n+1}z^{n+1}P_{n+1}(z)=f_{n+1}(z)-\a_i$ have the same number of zeros (counted with multiplicity) inside $B(0:r_{n+1})$. Note that $z=0$ is a zero of $f_n$ and $f_{n+1}$ of multiplicity 1 (since $f'_n(0)=f_{n+1}'(0)=1$). Suppose that $y_j$ ($j\in [1,s]$) is a zero of $f_n(z)-\a_i$ of multiplicity $m\geq 1$. Since $m \le \deg f_n$, $y_j$ is a zero of $f_{n+1}$ with multiplicity at least $m$. Thus, the sets $f_{n}^{-1}(\{\a_j\})\cap B(0:r_{n+1})$ and $f_{n+1}^{-1}(\{\a_j\})\cap B(0:r_{n+1})$ have the same cardinality and then they are equal, for $1\le j\le 3n+1$. This argument ensures that in our construction no new preimage under $f_{n+1}$ of $\{\a_1,\ldots, \a_{3n+1}\}$ lying in $B(0:r_{n+1})$ will appear apart from those ones under $f_n$. Note also that, since we will only choose $\e_{m+1}\in \A$, we have that, for every $n$, $f_{n+1}(\{\a_1,\ldots, \a_{3n+1}\})$ and $f_{n+1}^{-1}(\{\a_1,\ldots, \a_{3n+1}\})$ are subsets of $\QQ$. 

Now, we shall prove that it is possible to choose  a nonzero $\e_{n+1}\in \A$ satisfying (iv) and such that $a_{n+1}$ is a rational number (note that, by construction, the first $n$ coefficients of $f_{n+1}$ remains unchanged). In fact, let $c_{n+1}$ be the coefficient of $z^{n+1}$ in $f_n(z)$. We have that $a_{n+1}=c_{n+1}+P_{n+1}(0)\e_{n+1}$. Thus, one can choose a rational number $p/q$ such that
\[
0<|c_{n+1}-p/q|<\min\left\{\dfrac{|P_{n+1}(0)|}{L(P_{n+1})(n+1)^{n+1+\deg P_{n+1}}}, \dfrac{|P_{n+1}(0)|\displaystyle\min_{|z|=r_{n+1}}|f_n(z)-\a_i|}{\displaystyle\max_{|z|=r_{n+1}}|z^{n+1}P_{n+1}(z)|}\right\},
\]
for all $i\in [1,3n+1]$. Therefore, by defining $\e_{n+1}=(p/q-c_{n+1})/P_{n+1}(0)$, we get $a_{n+1}=p/q$ and $0<|\e_{n+1}|<\frac{1}{L(P_{n+1})(n+1)^{n+1+\deg P_{n+1}}}$ (also it satisfies (\ref{rouche0})).

Thus, by construction, the function $f(z)=z+\sum_{n\geq 2}\e_nz^nP_n(z)=\sum_{n\geq 1} a_nz^n$ is entire, $f(\QQ)\cup f^{-1}(\QQ)\subseteq \QQ$ and $a_n\in \Q$ as desired. Indeed, if $j\le 3n+1$, then $f_{n+1}(\a_j)=f_n(\a_j)$, so in particular $f_n(\a_j)=f_j(\a_j), \forall n\ge j$, and $f(\a_j)=\lim_{n\to \infty} f_n(\a_j)=f_j(\a_j)\in \QQ$. On the other hand, if $j\le 3n+1$ and $i\le n$ then $f_{n+1}^ {-1}(\a_j)\cap B(0,r_i)=f_n^ {-1}(\a_j)\cap B(0,r_i)$. Therefore, if $k=\max\{i,j\}$, we have $f_n^{-1}(\a_j)\cap B(0,r_i)=f_k^{-1}(\a_j)\cap B(0,r_i)$ for every $n\ge k$. This implies that $f^{-1}(\a_j)\cap B(0,r_i)\supseteq f_k^{-1}(\a_j)\cap B(0,r_i)$. Since $f$ is non-constant, we should have $f^{-1}(\a_j)\cap B(0,r_i)=f_k^{-1}(\a_j)\cap B(0,r_i)\subseteq \QQ$. Indeed, if there were another element $w$ of  $f^{-1}(\a_j)\cap B(0,r_i)$, it should be at positive distance of the finite set $f_k^{-1}(\a_j)\cap B(0,r_i)$, but, since $f=\lim_{n\to \infty} f_n$, arbitrarily close to $w$ there should be, for $n$ large, an element of $f_n^{-1}(\a_j)$ (again by Rouch\' e's theorem), which contradicts the equality $f_n^{-1}(\a_j)\cap B(0,r_i)=f_k^{-1}(\a_j)\cap B(0,r_i)$.

The proof that we can choose $f$ to be transcendental follows because there is an $\infty$-ary tree of different possibilities for $f$ (in each step we have infinitely many possible choices for $\e_{n+1}$, and so for $a_{n+1}$). Thus, we have constructed uncountably many possible functions, and the algebraic entire functions taking $\QQ$ into itself must be polynomials belonging to $\QQ[z]$, which is a countable subset.
\qed







\begin{thebibliography}{9999}


\bibitem{bookmahler} K. Mahler, \textit{Lectures on Transcendental Numbers}, Lecture Notes in Math., \textbf{546}, Springer-Verlag, Berlin, 1976.











\bibitem{19} P. St\"{a}ckel, Ueber arithmetische Eingenschaften analytischer Functionen, \textit{Math. Ann.} \textbf{46} (1895), 513--520.

\bibitem{20} P. St\"{a}ckel, Arithmetische eingenschaften analytischer Functionen, \textit{Acta Math.} \textbf{25} (1902), 371--383.


\bibitem{wal1} M. Waldschmidt, Algebraic values of analytic functions, Proceedings of the International Conference on Special Functions and their Applications (Chennai, 2002).  \textit{J. Comput. Appl. Math.} \textbf{160} (2003), 323--333.




\end{thebibliography}
\end{document}